\documentclass{article}
\usepackage{amsmath}
\usepackage{amssymb}
\usepackage{graphicx}
\begin{document}

\title{Schubert calculus and Intersection theory of Flag manifolds}
\author{Haibao Duan and Xuezhi Zhao}
\date{}
\maketitle

\begin{abstract}
Hilbert's 15th problem called for a rigorous foundation of Schubert's calculus, in which a long standing and challenging part is Schubert's problem of characteristics. In the course of securing the foundation of algebraic geometry, Van der Waerden and Andr\'{e} Weil attributed the problem to the determination of the intersection theory of flag manifolds.

This article surveys the background, content, and resolution of the problem of characteristics. Our main results are a unified formula for the characteristics, and a system description for the intersection rings of flag manifolds. We illustrate the effectiveness of the formula and the algorithm via explicit examples.

\bigskip

\noindent 2010 Mathematical Subject: Classification: 14M15; 57T15; 01A65

\end{abstract}

\section{Introduction}

Hilbert's 15th problem \cite{Hil} is an inspiring and far-reaching one. It promotes the enumerative geometry of the 19th century growing into the algebraic geometry founded by Van der Waerden and Andr\'{e} Weil, and makes Schubert calculus integrated deeply into many branches of mathematics. However, despite great many achievements in the 20th century (e.g. \cite{Ful,Kl3,So}), the part of the problem of finding an effective rule performing the calculus has been stagnant for a long time, notably, the Schubert's problem of characteristics \cite[\S 8]{VDW1}, or the Weil's problem \cite[p.331]{Weil} on the intersection theory of flag manifolds $G/P$, where $G$ is a compact connected Lie group and $P$ a parabolic subgroup \cite{BGG}.

In the series of works \cite{D2,D3,D5,DZ4,DZ6}, we have addressed both of the problem of characteristics and Weil's problem, implying that the 15th problem has been solved satisfactorily \cite[Remark 6.3]{DZ7}. The purpose of this article is to give an overview of the background, content and the resolution of Schubert's characteristics. In \S 2 we glimpse the development from Apollonius's work "Tangencies" to Lefschetz's homology theory, which reflects the evolution of the basic ideas from enumerative geometry to intersection theory. Section \S 3 summarises the pioneer contributions of Van der Waerden, Ehresmann, Weil, Bernstein-Gel'fand-Gel'fand to Schubert calculus, which led to a great clarification of Schubert's characteristics. Our main results are introduced in \S 4, where we present a formula expressing the characteristics of a flag manifold $G/P$ as a polynomial in the Cartan numbers of the Lie group $G$ (Theorem 4.5), develop a system description for the intersection rings of flag manifolds (Theorem 4.8), and illustrate the effectiveness of our programs in Examples 4.6, 4.9 and 4.12. In particular, since our approach uses the Cartan matrices of Lie groups as the main input, the solutions can be implemented successfully by computer programs, so that the intersection theory of flag manifolds becomes easily accessible by a broad readers.
\section{An introduction to intersection theory}

In the 2nd century B.C. Apollonius of Perga obtained the
following enumerative result in the paper ``Tangencies''.

\bigskip

\noindent \textbf{Apollonius's Theorem.} \textsl{The number of circles
tangent to three general circles in plane is $8$.}\hfill $\square$

\bigskip

\noindent The original proof of Apollonius was lost, but a record of the
theorem by Pappus dated in the 4th century survived. During the
Renaissance different proofs were founded respectively by Francois Viete, Adriaan van Roomen, Joseph
Diaz Gergonne and Isaac Newton \cite[p.159]{CB}. For a pictorial illustration of the theorem see the cover-page story of the book ``3264 and all that''\cite{EH}.

Descartes's discovery of the Euclidean coordinates makes it possible for
geometers (e.g. Maclaurin, Euler, Bezout) to exploit polynomial system to
characterize geometric figures that satisfy a system of incidence
conditions. Consequently, many enumerative problems admit the following
algebraic formulation.

\bigskip

\noindent \textbf{Problem 2.1.} \textsl{Given a system of $n$ polynomials in $n$ variables with complex coefficients}
\begin{enumerate}
\item[(2.1)] $\left\{
\begin{array}{c}
f_{1}(x_{1},\ldots,x_{n})=0 \\
\vdots \\
f_{n}(x_{1},\ldots,x_{n})=0%
\end{array}%
\right .$
\end{enumerate}
\noindent \textsl{find the number of solutions to the system.}
\hfill $\square$

\bigskip

Problem 2.1 is a fundamental one in algebra. In the case $n=1$ Gauss proved in the 1820's that the number of
zero's of a polynomial in a single variable is the degree of that
polynomial, well-known as the \textsl{fundamental theorem of algebra}.

Letting $g_{i}$ be the homogenization of the polynomial $f_{i}$ in (2.1), we get
in the $n$-dimensional complex projective space $\mathbb{C}P^{n}$ a
hypersurface $N_{i}:=g_{i}^{-1}(0)$. In general, the zero locus of a
homogeneous system on a complex projective space is called a \textsl{projective
variety}. Naturally, the study of problem 2.1 leads to the fundamental
problem of intersection theory.

\bigskip

\noindent \textbf{Problem 2.2.} \textsl{Given $k$ subvarieties $
N_{1},\ldots ,N_{k}$ in a smooth projective manifold $M$ that satisfy the
dimension constraint $\dim N_{1}+\cdots +\dim N_{k}=(k-1)\dim M$, find the
number $\left\vert N_{1}\cap \cdots \cap N_{k}\right\vert $ of the common
intersection points when the subvarieties $N_{i}$'s are in general position.}\hfill $\square$

%\textbf{Insert picture "intersection" here}

\begin{center}
\includegraphics[scale=0.4]{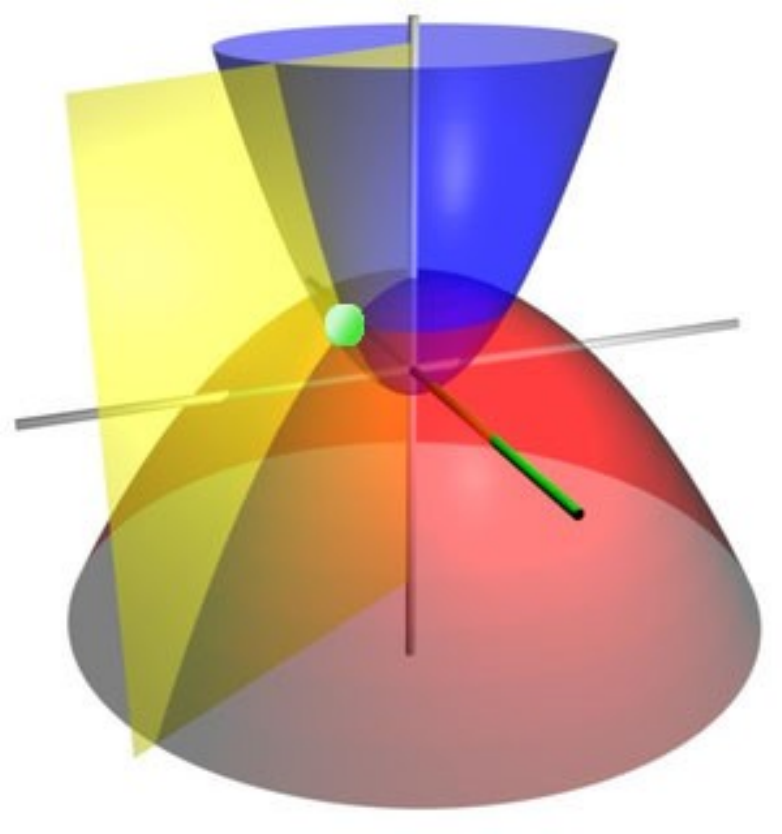}
\end{center}

In the course of studying Problem 2.2 S. Lefschetz developed the homology theory for the
cellular complexes \cite[1926]{Lef}. In the perspective of this theory let $\alpha_{i}\in H^{\dim M-\dim N_{i}}(M)$ be the Poincar\'{e} dual of the cycle class represented by the subvariety $N_{i}\subset M$.

\bigskip

\noindent \textbf{Problem 2.3. }\textsl{Given $k$ projective subvarieties $%
N_{1},\ldots ,N_{k}$ in a smooth projective manifold $M$ that satisfy the
dimension constraint $\dim N_{1}+\cdots +\dim N_{k}=(k-1)\dim M$, compute
the Kronecker pairing}

\begin{quote}
$\left\langle \alpha _{1}\cup \cdots \cup \alpha _{k},[M]\right\rangle =?$
\end{quote}

\noindent \textsl{where $\cup$ means the cup product on the cohomology ring $%
H^{\ast }(M)$, and where $[M]$ denotes the orientation class of $M$}. \hfill
$\square$

\bigskip

Through problems 2.1 to 2.3 we have briefly reviewed three seemingly different approaches to the problems of enumerative geometry. Given that the effective computability is the primary task of enumerative geometry, a natural question is:  which approach is the mostly calculable one? The development of the intersection theory shall tell us the answer.

\section{Schubert's problem of characteristics}
Hermann Schubert (1848-1911) received his Ph.D. from the University of
Halle, Germany in 1870. His doctoral thesis ``The theory of
characteristics'' \cite{Sch1} is about enumerative geometry. Prior to
this he had shown that there are $16$ spheres tangent to $4$ general spheres in space, a direct
extension of the Apollonius theorem.

In 1879 Schubert published the celebrated book ``Calculus of Enumerative
Geometry'' \cite{Sch2} that represents the summit of intersection theory in
the late 19th century \cite[p.2]{EH}. While developing M. Chasles's work on
conics \cite{Cha} he demonstrated amazing applications of intersection
theory to enumerative geometry, such as

\begin{quote}
i) The number of conics tangent to 8 general quadrics in space is 4,407,296;

ii) The number of quadrics tangent to 9 general quadrics in space is 666,841,088.

iii) The number of twisted cubic curves tangent to 12 general quadrics in space is
5,819,539,783,680.
\end{quote}

\noindent Nevertheless, in addition to the extensive use of the controversial ``principle of conservation of numbers'' \cite{Kl1,VDW}, Schubert's exposition was so sketching that gave ``no definition of
intersection multiplicity, no way to find it nor to calculate it'' \cite{Y}. At the outset of the 20th century Hilbert made finding rigorous foundations for Schubert calculus one of his celebrated problems, where he praised also the advantage of the calculus to foresee the final degree of a polynomial system before carrying out the actual process of elimination \cite{Hil}.

In order to gain insight into the central part of Schubert's approach to those spectacle enumerative numbers, we resort to the table of the characteristics for the variety of complete conics in space from his book \cite[p.95]{Sch2}:

\begin{center}
{\footnotesize
\begin{tabular}{l|l|l|l}
\hline
$\mu ^{3}\nu ^{5}=1$ & $\mu ^{2}\nu ^{6}=8$ & $\mu \nu ^{7}=34$ & $\nu
^{8}=92$ \\
$\mu ^{3}\nu ^{4}\rho =2$ & $\mu ^{2}\nu ^{5}\rho =14$ & $\mu \nu ^{6}\rho
=52$ & $\nu ^{7}\rho =116$ \\
$\mu ^{3}\nu ^{3}\rho ^{2}=4$ & $\mu ^{2}\nu ^{4}\rho ^{2}=24$ & $\mu \nu
^{5}\rho ^{2}=76$ & $\nu ^{6}\rho ^{2}=128$ \\
$\mu ^{3}\nu ^{2}\rho ^{3}=4$ & $\mu ^{2}\nu ^{3}\rho ^{3}=24$ & $\mu \nu
^{4}\rho ^{3}=72$ & $\nu ^{5}\rho ^{3}=104$ \\
$\mu ^{3}\nu \rho ^{4}=2$ & $\mu ^{2}\nu ^{2}\rho ^{4}=16$ & $\mu \nu
^{3}\rho ^{4}=48$ & $\nu ^{4}\rho ^{4}=64$ \\
$\mu ^{3}\rho ^{5}=1$ & $\mu ^{2}\nu \rho ^{5}=8$ & $\mu \nu ^{2}\rho
^{5}=24 $ & $\nu ^{3}\rho ^{5}=32$ \\
& $\mu ^{2}\rho ^{6}=4$ & $\mu \nu \rho ^{6}=12$ & $\nu ^{2}\rho ^{6}=16$ \\
&  & $\mu \rho ^{7}=6$ & $\nu \rho ^{7}=8$ \\
&  &  & $\rho ^{8}=4$ \\ \hline
\end{tabular}
}

\small {Table 1. The characteristics of the space of complete conics on $CP^{3}$.}
\end{center}

\noindent where the symbols $\mu, \nu$ and $\rho$ stand for the subvarieties of conics passing through a given point, intersecting a given line, and tangent to a given plane, respectively. The table consists of the equalities evaluating the monomials $\mu^{m}\nu ^{n}\rho ^{8-m-n}$ by integers, which were called \textsl{characteristics} by Schubert, and \textsl{the Schubert's symbolic
equations} by earlier researchers. Schubert emphasized that \textsl{the problem of characteristics} is the fundamental one of enumerative geometry \cite{Kl2,Sch1,Sch2,Sch4}. However, to state the problem in its natural simplicity and generality, one has to wait until 1950's for the celebrated ``basis theorem of Schubert calculus''. Let us recall the relevant works on the subject.

The study of the characteristics began with the Italian school headed by
Segre, Enriques and Severi. Two representing papers of the school are ``The
principle of conservation of numbers'' and ``The foundation of enumerative
geometry and the theory of characteristics'' due to Severi \cite{Sev,Sev1}.
Regarding these works Van de Wareden \cite{VDW} commented that they
``erected an admirable structure, but its logical foundation was shaky. The
notions were not well-defined, and the proofs were insufficient''.

In the pioneer work ``Topological foundation of enumerative geometry''
\cite[1930]{VDW1} Van der Waerden interpreted the Schubert's characteristics in the perspective of the homology theory developed by Lefschetz \cite{Lef} (e.g.Problem 2.3). He had the following crucial observations that enlightened the course of the later studies on the 15th problem:

\begin{quote}
1) Each Schubert's symbolic equation is a relation on the homology of a
projective manifold;

2) The solvability of Schubert's characteristic problem relies on a finite
basis of the homology of the relevant manifold;

3) The determination of the intersection products in homology is the goal of
all enumerative methods.
\end{quote}

\noindent C. Ehresmann \cite[1934]{Eh} went two important steps further: he discovered that

\begin{quote}
4) The parameter spaces of the geometric figures of Schubert are in
principle certain types of \textsl{flag manifolds} $G/P$ (see Remark 3.6);

5) For the Grassmannian $G_{n,k}$ of $k$-planes on the $n$-space $\mathbb{C}^{n}$ the Schubert symbols form exactly a basis of the homology $H_{\ast }(G_{n,k})$.
\end{quote}

In what follows we denote by $W(P,G)$ the set of left cosets of the Weyl group $W(G)$ of $G$ by the Weyl group $W(P)$ of $P$, and let $l:W(P,G)\rightarrow \mathbb{Z}$ be the associated length function \cite{BGG}. With the in-depth research on the structures of Lie groups (e.g. \cite {Bor1}) the vague term ``Schubert symbols'' in the early literature was gradually replaced by such rigorous defined geometric objects as ``Schubert cells'' or ``Schubert varieties''. In particular, extending Ehresmann's work \cite{Eh} on the Grassmannian $G_{n,k}$, the following result was announced by Chevalley \cite {Ch} for the complete flag manifold $G/T$ (where $T\subset G$ is a maximal torus), and extended to all flag manifolds $G/P$ by Bernstein-Gel'fand-Gel'fand \cite[Proposition 5.1]{BGG}.

\bigskip

\noindent \textbf{Theorem 3.1.} \textsl{Every flag manifold $G/P$ has a canonical
decomposition into the Schubert cells $S_{w}$, parameterized by the elements $w$ of $W(P,G)$,}

\begin{enumerate}
\item[(3.1)] $G/P=\underset{w\in W(P,G)}{\cup }S_{w}, \quad \dim S_{w}=2l(w),$
\end{enumerate}

\noindent \textsl{where the closure $X_{w}$ of each cell $S_{w}$ is a subvariety of $G/P$, called the Schubert variety on $G/P$ associated to $w\in W(P,G)$.}\hfill $\square$

\bigskip

Since only even dimensional cells are involved in the partition (3.1), the
set $\{[X_{w}], w\in W(P,G)\}$ of fundamental classes forms an additive
basis of the homology $H_{\ast }(G/P)$. The co-cycle classes $s_{w}\in
H^{\ast }(G/P)$ Kronecker dual to the basis (i.e. $\left\langle
s_{w},[X_{u}]\right\rangle =\delta _{w,u}, ~w,u\in W(P,G)$) gives rise to the
\textsl{Schubert class} associated to $w\in W(P,G)$. Theorem 3.1 implies the
following result, which is expected by Van der Waerden \cite[\S 8]{VDW1}, and is well-known as the ``basis theorem of Schubert calculus''.

\bigskip

\noindent \textbf{Theorem 3.2.}(\cite[Proposition 5.2]{BGG}) The set $\{s_{w}, w\in W(P,G)\}$ of Schubert
classes forms a basis of the cohomology $H^{\ast}(G/P)$. \hfill $\square$

\bigskip

An immediate consequence of the basis theorem is that any product $s_{u_{1}}\cdots s_{u_{k}}$ in the Schubert classes can be uniquely expressed as a linear combination of the basis elements

\begin{enumerate}
\item[(3.2)] $s_{u_{1}}\cdots s_{u_{k}}=\sum\limits_{w\in
W(P,G),l(w)=l(u_{1})+\cdots +l(w_{k}),\text{ }}c_{u_{1},\ldots
,u_{k}}^{w}\cdot s_{w}$, $c_{u_{1},\ldots ,u_{k}}^{w}\in \mathbb{Z}$
\end{enumerate}

\noindent by which Schubert's problem of characteristics \cite{Eh,Sch3,Sch4,VDW1} has the following concise expression \footnote{By Coolidge \cite[p.184]{Co} ``The fundamental problem which occupies Schubert is to express the product of two of these symbols in terms of others linearly. He succeeds in part.''}:

\bigskip

\noindent \textbf{Problem 3.3.} \textsl{Given any monomial $s_{u_{1}}\cdots
s_{u_{k}}$ in the Schubert classes, determine the characteristics numbers $c_{u_{1},\ldots ,u_{k}}^{w}$ in the linear expansion $(3.2)$}.\hfill $\square$

\bigskip

In the momentous treatise ``Foundations of Algebraic Geometry \cite{Weil}'' A. Weil completed the definition of intersection multiplicities for the first time, and summarized the task of Schubert calculus in the context of the modern intersection theory \footnote{``the classical Schubert calculus amounts to the determination of the intersection-rings on Grassmann varieties and on the so-called flag manifolds of projective geometry'' \cite[p.331]{Weil}, where we note that for a flag manifold $G/P$ the Chow ring $A^{\ast }(G/P)$ is canonically isomorphic to the cohomology $H^{\ast }(G/P)$}:

\bigskip

\noindent \textbf{Problem 3.4.} \textsl{Determine the intersection rings of flag manifolds $G/P$.} \hfill $\square$

\bigskip

Weil commented his problem as ``the modern form taken by the topic formerly
known as \textsl{enumerative geometry}'' \cite[p.331]{Weil}. We show that

\bigskip

\noindent \textbf{Theorem 3.5.} \textsl{For the flag manifolds $G/P$ the
Weil's problem is equivalent to the Schubert's one.}

\bigskip

\noindent \textbf{Proof.} A ring is an abelian group $R$ that is furnished
with a multiplication $R\times R\rightarrow R$. By the basis theorem the
cohomology $H^{\ast}(G/P)$ has a canonical basis consisting of Schubert
classes. Therefore, the multiplication on $H^{\ast}(G/P)$ is uniquely
determined by the product among the basis elements, which is handled by the
characteristics $c_{u_{1},\ldots ,u_{k}}^{w}$. \hfill $\square $

\bigskip

\noindent \textbf{Remark 3.6.} For the case $k=2$ the characteristics $%
c_{w_{1},\ldots ,w_{k}}^{w}$ admit various interpretations. They are called
\textsl{the Schubert's structure constants} of the flag manifold $G/P$ in topology; and
\textsl{the Littlewood-Richardson coefficients} in representation theory
\cite{LR}.

In certain cases the parameter spaces of the geometric figures concerned by
Schubert \cite[Chapt.IV]{Sch2} fail to be flag manifolds, but can be constructed by
performing finite steps of blow-ups on flag manifolds, see examples in Fulton \cite[Section 10.4]{Ful},
Eisenbud-Harris\cite[Chap.13]{EH}, or in \cite{DL} for the constructions of
the parameter spaces of the complete conics and quadrics on $\mathbb{C}P^{3}$. As results the relevant characteristics can be computed from those of flag manifolds via strict transformations (e.g
\cite[Examples 5.11; 5.12]{DL}). \hfill $\square $

\section{Intersection theory of flag manifolds}
To secure the foundation of a ``\textsl{calculus}'' it suffices to decide the objects to be calculated, and to determine accordingly the rules of the calculation (e.g. \cite[Chap.2]{B}, \cite{M}). As for the Schubert's calculus we have seen from \S 3 that the objects to be calculated have been clarified to be the Schubert symbols, or Schubert varieties. In this section we determine the rule of the calculus by a unified formula computing the characteristics, and apply the formula to complete the intersection theory of flag manifolds.
\subsection{Observation and expectation}

The major difficulties that one encounters with the problem of characteristics are fairly transparent:

\begin{quote}
i) The simply-connected simple Lie groups $G$ consist of the three infinite
families of classical Lie groups $Spin(n),Sp(n),SU(n)$, as well as the five
exceptional ones $G_{2},F_{4},E_{6},E_{7},E_{8}$;

ii) for a simple Lie group $G$ with rank $n$ there are precisely $2^{n}-1$
parabolic subgroups $P$ on $G$.
\end{quote}

\noindent That is, there exist plenty of flag manifolds $G/P$ whose
geometries and topologies vary considerably with respect to different
choices of $G$ and $P$. In addition

\begin{quote}
iii) the number of Schubert classes of $G/P$ agrees with the Euler
characteristic $\chi (G/P)$, which is normally very large,
\end{quote}

\noindent not to mention the number of the relevant characteristics. For
instance, for an exceptional Lie group $G$ with a maximal torus $T$ the
Euler characteristic $\chi (G/T)$ of the flag manifold $G/T$ is given in the table below

\begin{quote}
\begin{tabular}{l||l|l|l|l|l}
\hline
$G$ & $G_{2}$ & $F_{4}$ & $E_{6}$ & $E_{7}$ & $E_{8}$ \\ \hline
$\chi (G/T)$ & $12$ & $1152$ & $2^{7}\cdot 3^{4}\cdot 5$ & $2^{10}\cdot
3^{4}\cdot 5\cdot 7$ & $2^{14}\cdot 3^{5}\cdot 5^{2}\cdot 7$ \\ \hline
\end{tabular}%
.
\end{quote}

\noindent Summarizing, studies case by case can never reach a complete
solution to the problem.

On the other hand, according to E. Cartan's beautiful work on compact Lie groups, associate to each simple Lie group $G$ there is a Cartan matrix $C$, which acts the role of ``the cosmological constants''  to classify all flag manifolds $G/P$, see discussions in the coming section. As examples, for the five exceptional Lie groups those matrices are
{\renewcommand{\arraystretch}{1}
\setlength{\arraycolsep}{0.4 mm}
$$G_{2}:\left(
\begin{array}{cc}
2 & -1 \\
-3 & 2%
\end{array}%
\right), \ \ \
F_{4}:\left(
\begin{array}{cccc}
2 & -1 & 0 & 0 \\
-1 & 2 & -2 & 0 \\
0 & -1 & 2 & -1 \\
0 & 0 & -1 & 2%
\end{array}%
\right), \ \ \ E_{6}:\left(
\begin{array}{cccccc}
2 & 0 & -1 & 0 & 0 & 0 \\
0 & 2 & 0 & -1 & 0 & 0 \\
-1 & 0 & 2 & -1 & 0 & 0 \\
0 & -1 & -1 & 2 & -1 & 0 \\
0 & 0 & 0 & -1 & 2 & -1 \\
0 & 0 & 0 & 0 & -1 & 2%
\end{array}%
\right), $$
}
{\renewcommand{\arraystretch}{1}
\setlength{\arraycolsep}{0.4 mm}
$$E_{7}:\left(
\begin{array}{ccccccc}
2 & 0 & -1 & 0 & 0 & 0 & 0 \\
0 & 2 & 0 & -1 & 0 & 0 & 0 \\
-1 & 0 & 2 & -1 & 0 & 0 & 0 \\
0 & -1 & -1 & 2 & -1 & 0 & 0 \\
0 & 0 & 0 & -1 & 2 & -1 & 0 \\
0 & 0 & 0 & 0 & -1 & 2 & -1 \\
0 & 0 & 0 & 0 & 0 & -1 & 2%
\end{array}%
\right), \ \ \ E_{8}:\left(
\begin{array}{cccccccc}
2 & 0 & -1 & 0 & 0 & 0 & 0 & 0 \\
0 & 2 & 0 & -1 & 0 & 0 & 0 & 0 \\
-1 & 0 & 2 & -1 & 0 & 0 & 0 & 0 \\
0 & -1 & -1 & 2 & -1 & 0 & 0 & 0 \\
0 & 0 & 0 & -1 & 2 & -1 & 0 & 0 \\
0 & 0 & 0 & 0 & -1 & 2 & -1 & 0 \\
0 & 0 & 0 & 0 & 0 & -1 & 2 & -1 \\
0 & 0 & 0 & 0 & 0 & 0 & -1 & 2%
\end{array}%
\right). $$}
\noindent This raises the following question: can one express the characteristic numbers, as well as the intersection ring of a flag manifold $G/P$, merely in term of the Cartan matrix of the Lie group $G$?  In this section we fulfill this expectation.

\subsection{Numerical construction of a Weyl group}

Therefore, let $C=\left( c_{i,j}\right) _{n\times n}$ be the Cartan matrix
of some compact simple Lie group $G$, and let $\mathbb{R}^{n}$ be the $n-$dimensional real vector space with basis $\left\{ \omega _{1},\cdots
,\omega _{n}\right\}$. Define in terms of $C$ the endomorphisms $\sigma
_{i}\in End(\mathbb{R}^{n}),1\leq i\leq n,$ by the formula

\begin{quote}
$\sigma _{i}(\omega _{k})=\left\{
\begin{tabular}{l}
$\omega _{k}$ if $i\neq k$; \\
$\omega _{k}-(c_{k,1}\omega _{1}+c_{k,2}\omega _{2}+\cdots +c_{k,n}\omega
_{n})$ if $i\neq k$.%
\end{tabular}%
\right. $
\end{quote}

\noindent By general properties of Cartan matrices we have $\sigma _{i}^{2}=id$,
implying that $\sigma _{i}\in Aut(\mathbb{R}^{n})$. It can be further shown
that

\bigskip

\noindent \textbf{Lemma 4.1.} \textsl{The subgroup of $Aut(R^{n})$ generated
by the $\sigma _{i}$'s is isomorphic to the Weyl group $W(G)$ of $G$}.\hfill $\square$

\bigskip

For each subset $K\subset \{1,\ldots ,n\}$ there is a parabolic subgroup $%
P=P_{K}$, unique up to the conjugations on $G$, whose Weyl group $W(P)$ is
generated by those generators $\sigma _{{j}}\in W(G)$ with $j\notin K$. Resorting to the length
function $l$ on $W(G)$ we can furthermore embed the set $W(P;G)$ as the
subset of the group $W(G)$ \cite{BGG}

\begin{quote}
$W(P;G)=\{w\in W(G)\mid l(w_{{1}})\geq l(w)$, $w_{{1}}\in wW(P)\}$,
\end{quote}

\noindent and put $W^{m}(P;G):=\{w\in W(P;G)\mid l(w)=m\}$. By Lemma 4.1
each element $w\in W^{m}(P;G)$ admits a factorization of the form

\begin{quote}
$w=\sigma _{{i}_{{1}}}\circ \cdots \circ \sigma _{{i}_{{m}}}$ with $1\leq i_{%
{1}},\ldots ,i_{{m}}\leq n$,
\end{quote}

\noindent hence can be denoted by $w=\sigma _{I}$, where $I=(i_{{1}},\ldots
,i_{{m} })$. Such expressions of $w$ may not be unique, but the ambiguity
can be dispelled by employing the following notion. Furnish the set $
D(w):=\{I=(i_{{1}},\ldots,i_{{m}})\mid w=\sigma _{I}\}$ with the
lexicographical order $\preceq$ on the multi-indices $I$'s. We call a
decomposition $w=\sigma _{I}$ \textsl{minimized} if $I\in D(w)$ is the
minimal one. Clearly we have (e.g. \cite{Bu})

\bigskip

\noindent \textbf{Lemma 4.2.} \textsl{Every $w\in W(P;G)$ has a unique minimized decomposition.}\hfill $\square $

\bigskip

It follows that the set $W^{m}(P;G)$ is also ordered by the lexicographical
order $\preceq$ on the multi-index $I$'s, hence can be uniquely presented as

\begin{enumerate}
\item[(4.1)] $W^{m}(P;G)=\{w_{{m,i}}\mid 1\leq i\leq \beta (m)\}, ~
\beta (m):=\left\vert W^{m}(P;G)\right\vert $,
\end{enumerate}

\noindent where $w_{{m,i}}$ denotes the $i^{th}$ element in $W^{m}(P;G)$. In \cite{DZ1} the package ``\textsl{Decomposition}'' in MATHEMATICA is compiled, whose function is stated below.

\bigskip

\noindent \textbf{Algorithm I.} \textsl{Decomposition.}

\begin{quote}
\noindent \textbf{Input:} \textsl{The Cartan matrix $C=(a_{{ij}})_{{n\times n%
}}$ of $G$, and a subset $K\subset \{1,\ldots ,n\}$ to specify a parabolic
subgroup $P$.}

\textbf{Output:} \textsl{The set $W(P;G)$ being presented by the minimized
decompositions of its elements, together with the index system (4.1) imposed
by the order $\preceq$}.\hfill $\square$
\end{quote}

\noindent \textbf{Example 4.3.} Let $G=SU(n)$ be the special unitary group,
and let $k\in \{1,\ldots ,n-1\}$. The flag manifold $G/P_{\{k\}}$ is the
Grassmannian manifold $G_{n,k}$. Applying the \textsl{Decomposition} to the case $G_{9,4}$ we obtain the following table of minimized decompositions, as well as the order imposed by (4.1), for the elements $w\in W(P_{4};SU(9))$ with $l(w)\leq 8$.

\bigskip
 %% all cells
{\setlength{\tabcolsep}{0.6 mm}
\begin{tabular}{|l|l||l|l||l|l|}
\hline
$w_{i,j}$ & decomposition & $w_{i,j}$ & decomposition  & $w_{i,j}$ & decomposition \\
\hline

$w_{1,1}$  &
$[{4}]$ &

$w_{2,1}$  &
$[{3, 4}]$ &

$w_{2,2}$  &
$[{5, 4}]$\\

$w_{3,1}$  &
$[{2, 3, 4}]$ &

$w_{3,2}$  &
$[{3, 5, 4}]$ &

$w_{3,3}$  &
$[{6, 5, 4}]$\\

$w_{4,1}$  &
$[{1, 2, 3, 4}]$ &

$w_{4,2}$  &
$[{2, 3, 5, 4}]$ &

$w_{4,3}$  &
$[{3, 6, 5, 4}]$\\

$w_{4,4}$  &
$[{4, 3, 5, 4}]$ &

$w_{4,5}$  &
$[{7, 6, 5, 4}]$ &

$w_{5,1}$  &
$[{1, 2, 3, 5, 4}]$\\

$w_{5,2}$  &
$[{2, 3, 6, 5, 4}]$ &

$w_{5,3}$  &
$[{2, 4, 3, 5, 4}]$ &

$w_{5,4}$  &
$[{3, 7, 6, 5, 4}]$\\

$w_{5,5}$  &
$[{4, 3, 6, 5, 4}]$ &

$w_{5,6}$  &
$[{8, 7, 6, 5, 4}]$ &

$w_{6,1}$  &
$[{1, 2, 3, 6, 5, 4}]$\\

$w_{6,2}$  &
$[{1, 2, 4, 3, 5, 4}]$ &

$w_{6,3}$  &
$[{2, 3, 7, 6, 5, 4}]$ &

$w_{6,4}$  &
$[{2, 4, 3, 6, 5, 4}]$\\

$w_{6,5}$  &
$[{3, 2, 4, 3, 5, 4}]$ &

$w_{6,6}$  &
$[{3, 8, 7, 6, 5, 4}]$ &

$w_{6,7}$  &
$[{4, 3, 7, 6, 5, 4}]$\\

$w_{6,8}$  &
$[{5, 4, 3, 6, 5, 4}]$ &

$w_{7,1}$  &
$[{1, 2, 3, 7, 6, 5, 4}]$ &

$w_{7,2}$  &
$[{1, 2, 4, 3, 6, 5, 4}]$\\

$w_{7,3}$  &
$[{1, 3, 2, 4, 3, 5, 4}]$ &

$w_{7,4}$  &
$[{2, 3, 8, 7, 6, 5, 4}]$ &

$w_{7,5}$  &
$[{2, 4, 3, 7, 6, 5, 4}]$\\

$w_{7,6}$  &
$[{2, 5, 4, 3, 6, 5, 4}]$ &

$w_{7,7}$  &
$[{3, 2, 4, 3, 6, 5, 4}]$ &

$w_{7,8}$  &
$[{4, 3, 8, 7, 6, 5, 4}]$\\

$w_{7,9}$  &
$[{5, 4, 3, 7, 6, 5, 4}]$ &

$w_{8,1}$  &
$[{1, 2, 3, 8, 7, 6, 5, 4}]$ &

$w_{8,2}$  &
$[{1, 2, 4, 3, 7, 6, 5, 4}]$\\

$w_{8,3}$  &
$[{1, 2, 5, 4, 3, 6, 5, 4}]$ &

$w_{8,4}$  &
$[{1, 3, 2, 4, 3, 6, 5, 4}]$ &

$w_{8,5}$  &
$[{2, 1, 3, 2, 4, 3, 5, 4}]$\\

$w_{8,6}$  &
$[{2, 4, 3, 8, 7, 6, 5, 4}]$ &

$w_{8,7}$  &
$[{2, 5, 4, 3, 7, 6, 5, 4}]$ &

$w_{8,8}$  &
$[{3, 2, 4, 3, 7, 6, 5, 4}]$\\

$w_{8,9}$  &
$[{3, 2, 5, 4, 3, 6, 5, 4}]$ &

$w_{8,10}$  &
$[{5, 4, 3, 8, 7, 6, 5, 4}]$ &

$w_{8,11}$  &
$[{6, 5, 4, 3, 7, 6, 5, 4}]$\\

\hline
\end{tabular}
}
\bigskip

\noindent For more examples of the results produced by \textsl{Decomposition}
we refer to \cite[Sections 1.1--7.1]{DZ5}.

Geometrically, for any $w\in W(P;G)$ the Schubert variety $X_{w}$ can be explicitly constructed in terms of its minimized decomposition \cite{BGG,D3}. \hfill $\square$

\subsection{\noindent A unified formula for Schubert's characteristics}

Given an element $w\in W(P_{K};G)$ with minimized decomposition

\begin{quote}
$w=\sigma _{i_{1}}\circ \sigma _{i_{2}}\circ \cdots \circ \sigma _{i_{m}}$, $%
1\leq i_{1},i_{2},\cdots ,i_{m}\leq n$,
\end{quote}

\noindent \textsl{the structure matrix} of $w$ is the strictly
upper triangular matrix $A_{w}=\left( a_{s,t}\right) _{m\times m}$ defined
by the Cartan matrix $C=\left( c_{i,j}\right) _{n\times n}$ of $G$ as

\begin{quote}
$a_{s,t}=0$ if $s\geq t$, $-c_{i_{s},i_{t}}$ if $s<t$.
\end{quote}

\noindent As examples recall that the Cartan matrix of the exceptional Lie
group $G_{2}$ is

\begin{quote}
$C=\left(
\begin{array}{cc}
2 & -1 \\
-3 & 2%
\end{array}%
\right) $.
\end{quote}

\noindent By Lemma 4.1 the Weyl group $W(G_{2})$ has two generators $
\sigma_{1},\sigma_{2}$. Consider the following elements with length $4$:

\begin{quote}
$u=\sigma_{1}\circ \sigma_{2}\circ \sigma_{1}\circ \sigma_{2}$ and $%
v=\sigma_{2}\circ \sigma_{1}\circ \sigma_{2}\circ \sigma_{1}$.
\end{quote}

\noindent From the Cartan matrix $C$ one reads

\begin{quote}
$A_{u}=\left(
\begin{array}{cccc}
0 & 1 & -2 & 1 \\
0 & 0 & 3 & -2 \\
0 & 0 & 0 & 1 \\
0 & 0 & 0 & 0%
\end{array}%
\right) $ and $A_{v}=\left(
\begin{array}{cccc}
0 & 3 & -2 & 3 \\
0 & 0 & 1 & -2 \\
0 & 0 & 0 & 3 \\
0 & 0 & 0 & 0%
\end{array}%
\right) $.
\end{quote}

Let $\mathbb{Z}[x_{{1}},\ldots ,x_{{m}}]$ be the ring of polynomials in $x_{{%
1}},\ldots ,x_{{m}}$ that is graded by $\deg x_{i}=1$, and let $\mathbb{Z}%
[x_{{1}},\ldots ,x_{{m}}]^{(m)}$ be its subgroup spanned by the monomials with degree $m$. Given a $m\times m$ strictly upper triangular integer matrix $A=(a_{{i,j}})$ \textsl{the triangular operator $T_{{A}}$}
associated to $A$ is the linear map

\begin{quote}
$T_{A}:\mathbb{Z}[x_{{1}},\ldots ,x_{{m}}]^{(m)}\rightarrow \mathbb{Z}$
\end{quote}

\noindent defined recursively by the following elimination rules:

\begin{quote}
i) If $m=1$ (i.e. $A=(0)$) then $T_{{A}}(x_{{1}})=1$;

ii) If $h\in \mathbb{Z}[x_{{1}},\ldots ,x_{{m-1}}]^{(m)}$ then $T_{{A}}(h)=0$;

iii) If $h\in \mathbb{Z}[x_{{1}},\ldots ,x_{{m-1}}]^{(m-r)}$ with $r\geq 1$
then

$\qquad T_{{A}}(h\cdot x_{{m}}^{r})=T_{{A}_{1}}(h\cdot (a_{{1,m}}x_{{1}%
}+\cdots +a_{{m-1,m}}x_{{m-1}})^{r-1}),$
\end{quote}

\noindent where $A_{1}$ is the $(m-1)\times (m-1)$ (strictly upper
triangular matrix) obtained from $A$ by deleting both of the $m^{th}$ column
and row. Since every polynomial $h\in \mathbb{Z}[x_{{1}},\ldots ,x_{{m}%
}]^{(m)}$ admits the unique expansion

\begin{quote}
$h=\sum\limits_{{0\leq r\leq m}}h_{{r}}\cdot x_{{m}}^{r}$ with $h_{{r}}\in
\mathbb{Z}[x_{{1}},\ldots,x_{{m-1}}]^{(m-r)}$,
\end{quote}

\noindent the operator $T_{{A}}$ is well-defined by the rules i), ii) and
iii). It follows that

\bigskip

\noindent \textbf{Lemma 4.4.} \textsl{For any polynomial $h\in \mathbb{Z}[x_{{1}%
},\ldots ,x_{{m}}]^{(m)}$, the number $T_{A}(h)$ is a polynomial in the
entries of the matrix $A$ with degree $m$}. \hfill $\square $

\bigskip

Extending the main results of \cite{D2,D3,D5} we have shown in \cite[Theorem 2.4]{DZ7}
the following formula that expresses the Schubert's characteristics of a flag manifold $G/P$ as polynomials in the Cartan numbers of the group $G$.

\bigskip

\noindent \textbf{Theorem 4.5.} \textsl{Let $w\in W(P;G)$ be an element with
minimized decomposition $\sigma _{{i}_{{1}}}\circ \cdots \circ \sigma _{{i}_{%
{m}}}$ and structure matrix $A_{{w}}$. For any monomial $s_{u_{1}}\cdots s_{u_{k}}$ in the Schubert classes
with total degree $m$ one has}

\begin{enumerate}
\item[(4.2)] $c_{u_{1},\ldots ,u_{k}}^{w}=T_{A_{{w}}}\left( \underset{%
i=1,\ldots ,k}{\prod}\left( \underset{\sigma _{{I}}=u_{i},\left\vert
I\right\vert =l(u_{i}),I\subseteq \{1,\ldots ,m\}}{\sum }x_{I}\right)
\right) $\textsl{,}
\end{enumerate}

\noindent \textsl{where for a multi-index $I=\{j_{{1}},\ldots ,j_{{t}}\}$ we have set $\left\vert I\right\vert :=t$ and}

\begin{quote}
$\sigma _{{I}}:=\sigma _{{i}_{{j}_{{1}}}}\circ \cdots \circ \sigma _{{i}_{{j}%
_{{t}}}}\in W(G),~~x_{I}:=x_{i_{j_{1}}}\cdots x_{i_{j_{t}}}\in \mathbb{Z}[x_{%
{1}},\ldots,x_{{m}}]$.\hfill $\square$
\end{quote}

Since the matrix $A_{{w}}$ is constructed from the Cartan matrix of the group $G$ in term of the minimized decomposition of $w$, while the operator $T_{A_{{w}}}$ is evaluated easily by the elimination rules i)-iii) stated above, the formula (4.2) indicates an effective algorithm to evaluate the numbers $c_{u_{1},\ldots ,u_{k}}^{w}$. Combining these ideas the package ``\textsl{Characteristics}'' in MATHEMATICA has been compiled (e.g. \cite{DZ1}) whose function is described as follows.

\bigskip

\noindent \textbf{Algorithm II: }\textsl{Characteristics.}

\begin{quote}
\textbf{Input:} \textsl{The Cartan matrix $C=(a_{{ij}})_{{n\times n}}$ of $G$, and a subset $K\subset \{1,\ldots ,n\}$ to specify a parabolic subgroup $P$}.

\textbf{Output: }\textsl{The characteristics $c_{u_{1},\ldots ,u_{k}}^{w}$
of $G/P$.} \hfill $\square $
\end{quote}

\noindent \textbf{Example 4.6. The characteristics of the Schubert monomials at the top degree.} Let $G/P$ be a flag manifold with $\dim_{\mathbb C}G/P=m$. According to the basis theorem there exists a unique element $w_{0}\in W(P;G)$ so that $l(w_{0})=m$, and that the Schubert class $s_{w_{0}}$ generates the top degree cohomology $H^{2m}(G/P)=\mathbb{Z}$. It follows that, for any monomial $s_{u_{1}}\cdots s_{u_{k}}$ in the Schubert classes with total degree $m$, the characteristic number $c_{u_{1},\ldots ,u_{k}}^{w}$ is given by the equality (see in Problem 2.3):
$$\langle s_{u_{1}}\cdots s_{u_{k}},[G/P]\rangle =c_{u_{1},\ldots ,u_{k}}^{w_{0}},$$
\noindent which will be abbreviated by $s_{u_{1}}\cdots s_{u_{k}}=c_{u_{1},\ldots ,u_{k}}^{w_{0}}.$ In addition, for an element $w\in W(P;G)$ with minimized decomposition $w=\sigma_{I}$ we can use the notion $s_{I}$ to denote the Schubert class $s_{w}$.

The cohomology of the Grassmannians $G_{n,k}$ are the most classical and archetypal examples of intersection theory \cite{Sch3,Sch4}, \cite[p.4]{EH}. Traditionally, the characteristics $c_{u_{1},u_{2}}^{w}$ are given by the combinatorial \textsl{Littlewood-Richardson rule} \cite{LR}, rather than a closed formula. In contrast our formula (4.2) is practical for numerical computation. In term of the convention above we set
$$c_{r}:=s_{\{k-r+1,k-r+2,\cdots,k\}}\in H^{2r}(G_{n,k}), r=1,\cdots,k.$$
Then $c_{r}$ is also the $r^{th}$ Chern class of the canonical $k$-dimensional complex vector bundles on $G_{n,k}$ \cite{MS}. Applying \textsl{Characteristics} to the case $G_{9,4}$ we obtain the following table of characteristics for the monomials in the Chern classes at the top degree $\dim_{\mathbb C} G_{9,4}=20$ .

\bigskip
{\footnotesize  \renewcommand{\arraystretch}{0.8}

\begin{tabular}{|l|l|l|l|l|l||l|l|l|l|}
  \hline

$c_4^{5} =  1$   &

$c_3^{4}c_4^{2} =  1$   &

$c_2 c_3^{2}c_4^{3} =  1$   &

$c_2 c_3^{6} =  9$ \\ \hline

$c_2^{2}c_4^{4} =  1$   &

$c_2^{2}c_3^{4}c_4  =  6$   &

$c_2^{3}c_3^{2}c_4^{2} =  4$   &

$c_2^{4}c_4^{3} =  3$ \\ \hline

$c_2^{4}c_3^{4} =  45$   &

$c_2^{5}c_3^{2}c_4  =  26$   &

$c_2^{6}c_4^{2} =  16$   &

$c_2^{7}c_3^{2} =  231$ \\ \hline

$c_2^{8}c_4  =  126$   &

$c_2^{10} =  1296$   &

$c_1 c_3 c_4^{4} =  1$   &

$c_1 c_3^{5}c_4  =  4$ \\ \hline

$c_1 c_2 c_3^{3}c_4^{2} =  3$   &

$c_1 c_2^{2}c_3 c_4^{3} =  2$   &

$c_1 c_2^{2}c_3^{5} =  29$   &

$c_1 c_2^{3}c_3^{3}c_4  =  17$ \\ \hline

$c_1 c_2^{4}c_3 c_4^{2} =  10$   &

$c_1 c_2^{5}c_3^{3} =  141$   &

$c_1 c_2^{6}c_3 c_4  =  76$   &

$c_1 c_2^{8}c_3  =  756$ \\ \hline

$c_1^{2}c_3^{2}c_4^{3} =  2$   &

$c_1^{2}c_3^{6} =  19$   &

$c_1^{2}c_2 c_4^{4} =  1$   &

$c_1^{2}c_2 c_3^{4}c_4  =  12$ \\ \hline

$c_1^{2}c_2^{2}c_3^{2}c_4^{2} =  7$   &

$c_1^{2}c_2^{3}c_4^{3} =  4$   &

$c_1^{2}c_2^{3}c_3^{4} =  89$   &

$c_1^{2}c_2^{4}c_3^{2}c_4  =  48$ \\ \hline

$c_1^{2}c_2^{5}c_4^{2} =  26$   &

$c_1^{2}c_2^{6}c_3^{2} =  451$   &

$c_1^{2}c_2^{7}c_4  =  231$   &

$c_1^{2}c_2^{9} =  2556$ \\ \hline

$c_1^{3}c_3^{3}c_4^{2} =  6$   &

$c_1^{3}c_2 c_3 c_4^{3} =  3$   &

$c_1^{3}c_2 c_3^{5} =  59$   &

$c_1^{3}c_2^{2}c_3^{3}c_4  =  32$ \\ \hline

$c_1^{3}c_2^{3}c_3 c_4^{2} =  17$   &

$c_1^{3}c_2^{4}c_3^{3} =  276$   &

$c_1^{3}c_2^{5}c_3 c_4  =  141$   &

$c_1^{3}c_2^{7}c_3  =  1491$ \\ \hline

$c_1^{4}c_4^{4} =  1$   &

$c_1^{4}c_3^{4}c_4  =  24$   &

$c_1^{4}c_2 c_3^{2}c_4^{2} =  12$   &

$c_1^{4}c_2^{2}c_4^{3} =  6$ \\ \hline

$c_1^{4}c_2^{2}c_3^{4} =  175$   &

$c_1^{4}c_2^{3}c_3^{2}c_4  =  89$   &

$c_1^{4}c_2^{4}c_4^{2} =  45$   &

$c_1^{4}c_2^{5}c_3^{2} =  886$ \\ \hline

$c_1^{4}c_2^{6}c_4  =  436$   &

$c_1^{4}c_2^{8} =  5112$   &

$c_1^{5}c_3 c_4^{3} =  4$   &

$c_1^{5}c_3^{5} =  119$ \\ \hline

$c_1^{5}c_2 c_3^{3}c_4  =  59$   &

$c_1^{5}c_2^{2}c_3 c_4^{2} =  29$   &

$c_1^{5}c_2^{3}c_3^{3} =  539$   &

$c_1^{5}c_2^{4}c_3 c_4  =  264$ \\ \hline

$c_1^{5}c_2^{6}c_3  =  2962$   &

$c_1^{6}c_3^{2}c_4^{2} =  19$   &

$c_1^{6}c_2 c_4^{3} =  9$   &

$c_1^{6}c_2 c_3^{4} =  339$ \\ \hline

$c_1^{6}c_2^{2}c_3^{2}c_4  =  164$   &

$c_1^{6}c_2^{3}c_4^{2} =  79$   &

$c_1^{6}c_2^{4}c_3^{2} =  1744$   &

$c_1^{6}c_2^{5}c_4  =  832$ \\ \hline

$c_1^{6}c_2^{7} =  10302$   &

$c_1^{7}c_3^{3}c_4  =  104$   &

$c_1^{7}c_2 c_3 c_4^{2} =  49$   &

$c_1^{7}c_2^{2}c_3^{3} =  1047$ \\ \hline

$c_1^{7}c_2^{3}c_3 c_4  =  496$   &

$c_1^{7}c_2^{5}c_3  =  5912$   &

$c_1^{8}c_4^{3} =  14$   &

$c_1^{8}c_3^{4} =  641$ \\ \hline

$c_1^{8}c_2 c_3^{2}c_4  =  300$   &

$c_1^{8}c_2^{2}c_4^{2} =  140$   &

$c_1^{8}c_2^{3}c_3^{2} =  3437$   &

$c_1^{8}c_2^{4}c_4  =  1600$ \\ \hline

$c_1^{8}c_2^{6} =  20887$   &

$c_1^{9}c_3 c_4^{2} =  84$   &

$c_1^{9}c_2 c_3^{3} =  2025$   &

$c_1^{9}c_2^{2}c_3 c_4  =  936$ \\ \hline

$c_1^{9}c_2^{4}c_3  =  11853$   &

$c_1^{10}c_3^{2}c_4  =  552$   &

$c_1^{10}c_2 c_4^{2} =  252$   &

$c_1^{10}c_2^{2}c_3^{2} =  6792$ \\ \hline

$c_1^{10}c_2^{3}c_4  =  3102$   &

$c_1^{10}c_2^{5} =  42597$   &

$c_1^{11}c_3^{3} =  3927$   &

$c_1^{11}c_2 c_3 c_4  =  1782$ \\ \hline

$c_1^{11}c_2^{3}c_3  =  23892$   &

$c_1^{12}c_4^{2} =  462$   &

$c_1^{12}c_2 c_3^{2} =  13497$   &

$c_1^{12}c_2^{2}c_4  =  6072$ \\ \hline

$c_1^{12}c_2^{4} =  87417$   &

$c_1^{13}c_3 c_4  =  3432$   &

$c_1^{13}c_2^{2}c_3  =  48477$   &

$c_1^{14}c_3^{2} =  27027$ \\ \hline

$c_1^{14}c_2 c_4  =  12012$   &

$c_1^{14}c_2^{3} =  180609$   &

$c_1^{15}c_2 c_3  =  99099$   &

$c_1^{16}c_4  =  24024$ \\ \hline

$c_1^{16}c_2^{2} =  375804$   &

$c_1^{17}c_3  =  204204$   &

$c_1^{18}c_2  =  787644$   &

$c_1^{20} =  1662804$ \\ \hline

\end{tabular}
}
\begin{center}
\small{Table 2. The characteristics of the Grassmaniann $G_{9,4}$}
\end{center}

The \textsl{Characteristics} works equally well for other types of flag manifolds. For example consider the flag manifold $E_{6}/P_{\{2\}}$, where $P_{\{2\}}=S^{1}\cdot SU(6)$. Following Bourbaki's numbering of simple roots \cite{Bu} let $y_{1},y_{3},y_{4},y_{6}$ be respectively the Schubert classes $s_{I}$ with
$$I=\{2\}, \{5,4,2\}, \{6,5,4,2\}, \{1,3,6,5,4,2\},$$
Then the cohomology $H^{*}(E_{6}/P_{\{2\}})$ is generated by $y_{1},y_{3},y_{4},y_{6}$ by \cite[Theorem 3]{DZ4}. Applying \textsl{Characteristics} we obtain the following table of characteristics for all the monomials in the Schubert generators $y_{i}'s$ at the top degree $\dim_{\mathbb C} E_{6}/P_{\{2\}}=21$.

\begin{center}
{\footnotesize  \renewcommand{\arraystretch}{0.8}
\begin{tabular}{|l|l|l|l|l|l||l|l|l|l|}
 \hline

$y_3 y_6^{3} =  3$   &

$y_3 y_4^{3}y_6  =  3$   &

$y_3^{3}y_6^{2} =  21$   &

$y_3^{3}y_4^{3} =  21$  &

$y_3^{5}y_6  =  156$   \\ \hline

$y_3^{7} =  1158$   &

$y_1 y_4^{2}y_6^{2} =  2$   &

$y_1 y_4^{5} =  2$ &

$y_1 y_3^{2}y_4^{2}y_6  =  14$   &

$y_1 y_3^{4}y_4^{2} =  100$   \\ \hline

$y_1^{2}y_3 y_4 y_6^{2} =  9$   &

$y_1^{2}y_3 y_4^{4} =  9$ &

$y_1^{2}y_3^{3}y_4 y_6  =  66$   &

$y_1^{2}y_3^{5}y_4  =  483$   &

$y_1^{3}y_6^{3} =  6$   \\ \hline

$y_1^{3}y_4^{3}y_6  =  6$ &

$y_1^{3}y_3^{2}y_6^{2} =  42$   &

$y_1^{3}y_3^{2}y_4^{3} =  42$   &

$y_1^{3}y_3^{4}y_6  =  312$   &

$y_1^{3}y_3^{6} =  2328$ \\ \hline

$y_1^{4}y_3 y_4^{2}y_6  =  28$   &

$y_1^{4}y_3^{3}y_4^{2} =  201$   &

$y_1^{5}y_4 y_6^{2} =  18$   &

$y_1^{5}y_4^{4} =  18$ &

$y_1^{5}y_3^{2}y_4 y_6  =  132$   \\ \hline

$y_1^{5}y_3^{4}y_4  =  972$   &

$y_1^{6}y_3 y_6^{2} =  84$   &

$y_1^{6}y_3 y_4^{3} =  84$ &

$y_1^{6}y_3^{3}y_6  =  624$   &

$y_1^{6}y_3^{5} =  4677$   \\ \hline

$y_1^{7}y_4^{2}y_6  =  56$   &

$y_1^{7}y_3^{2}y_4^{2} =  404$ &

$y_1^{8}y_3 y_4 y_6  =  264$   &

$y_1^{8}y_3^{3}y_4  =  1956$   &

$y_1^{9}y_6^{2} =  168$   \\ \hline

$y_1^{9}y_4^{3} =  168$ &

$y_1^{9}y_3^{2}y_6  =  1248$   &

$y_1^{9}y_3^{4} =  9390$   &

$y_1^{10}y_3 y_4^{2} =  813$   &

$y_1^{11}y_4 y_6  =  528$ \\ \hline

$y_1^{11}y_3^{2}y_4  =  3936$   &

$y_1^{12}y_3 y_6  =  2496$   &

$y_1^{12}y_3^{3} =  18837$   &

$y_1^{13}y_4^{2} =  1638$ &

$y_1^{14}y_3 y_4  =  7917$   \\ \hline

$y_1^{15}y_6  =  4992$   &

$y_1^{15}y_3^{2} =  37752$   &

$y_1^{17}y_4  =  15912$ &

$y_1^{18}y_3  =  75582$   &

$y_1^{21} =  151164$  \\ \hline
\end{tabular}
}
\small {Table 3. The characteristics of the flag manifold $E_{6}/S^{1}\cdot SU(6)$.}
\end{center}

The contents of Tables 2 and 3 are compatible with Schubert's computation in Table 1. Let $M$ be the variety of complete conics in 3-space. Then $\dim_{\mathbb C} M=8$, while the ring $H^{*}(M)$ is generated by the Schubert's symbols $\mu,\rho,\nu\in H^{2}(M)$ \cite{DL}. That is, the equalities in Table 1 consist of the characteristics of the monomials $\mu^{m}\nu ^{n}\rho ^{8-m-n}$ in the symbols $\mu,\rho$ and $\nu$ at the top degree $8$. \hfill $\square$

\bigskip

\noindent \textbf{Remark 4.7.} Geometrically, the operator $T_{A_{{w}}}$ in (4.2) handles the integrations along the Schubert cell $X_{w}$ \cite{BGG,D3}.

The formula (4.2) has been extended in \cite{D4} to compute the products of the basis elements of the Grothendieck $K$-theory of flag manifolds, and to evaluate the Steenrod operations on Schubert classes in \cite{DZ2} .\hfill $\square$
\subsection{The intersection rings of flag manifolds}

As in Example 4.6 let $c_{i}\in H^{2i}(G_{{n,k}})$ be the $i^{th}$-Chern class. Borel \cite{Bor} has shown that

\begin{enumerate}
\item[(4.3)] $H^{\ast }(G_{{n,k}})=\mathbb{Z}[c_{{1}},\ldots ,c_{{k}%
}]/\left\langle c_{{n-k+1}}^{-1},\ldots,c_{{n}}^{-1}\right\rangle $,
\end{enumerate}

\noindent where $c_{{j}}^{-1}$ denotes the component of the formal inverse
of $1+c_{{1}}+\cdots +c_{{k}}$ in degree $j$, and where $\langle\cdots\rangle$ denotes the ideal generated by the enclosed polynomials. Comparing formula (4.3) with the contents in Table 2 reveals the following phenomena: the characteristic numbers are essential for enumerative geometry, but fails to be a concise way to characterize the structure of the ring $H^{\ast }(G_{{n,k}})$. It is the Weil's problem that motivates us the following extension of Borel's formula (4.3) to all flag manifolds.

\bigskip

\noindent \textbf{Theorem 4.8.} \textsl{For each flag manifold $G/P$ there
exist Schubert classes $y_{1},\cdots ,y_{n}$ such that}

\begin{enumerate}
\item[(4.4)] $H^{\ast }(G/P)=\mathbb{Z}[y_{1},\cdots ,y_{n}]/\left\langle
f_{1},\cdots ,f_{m}\right\rangle ,$
\end{enumerate}

\noindent \textsl{where $f_{i}\in \mathbb{Z}[y_{1},\cdots ,y_{n}], 1\leq i\leq m$, and where the numbers $n$ and $m$ are minimum subject to the presentation.}.

\bigskip

\noindent \textbf{Proof.} Let $H^{+}(G/P)$ be the subring of the cohomology $H^{*}(G/P)$ spanned by the homogeneous elements of positive degrees, and let $D(H^{*}(G/P))$ be the ideal of the decomposable elements of the ring $H^{+}(G/P)$. Since the cohomology $H^{*}(G/P)$ is torsion free and has a basis consisting of Schubert classes, there exist Schubert classes $y_{1},\cdots,y_{n}$ on $G/P$ that correspond to a basis of the quotient group $H^{+}(G/P)/D(H^{*}(G/P))$. It follows that the inclusion $\left\{ y_{1},\cdots ,y_{n}\right\} \subset H^{\ast }(G/P)$ induces a ring epimorphism

$$\pi:\mathbb{Z}[y_{1},\cdots,y_{n}]\rightarrow H^{\ast }(G/P).$$

\noindent By the Hilbert's basis theorem there exist finitely many polynomials $f_{1},\cdots ,f_{m}\in\mathbb{Z}[y_{1},\cdots,y_{n}]$ such that $\ker \pi=\left\langle f_{1},\cdots ,f_{m}\right\rangle$. We can of course assume that the number $m$ is minimum with respect to the formula (4.4).

As the cardinality of a basis of the quotient group $H^{+}(G/P)/D(H^{*}(G/P))$ the number $n$ is an invariant of $G/P$. In addition, if one changes the generators $y_{1},\cdots ,y_{n}$ to $y_{1}',\cdots ,y_{n}'$, then each old generator $y_{i}$ can be expressed as a polynomial $g_{i}$ in the new ones $y_{1}',\cdots,y_{n}'$. The invariance of the number $m$ is shown by the presentation

$$H^{\ast }(G/P)=\mathbb{Z}[y_{1}',\cdots ,y_{n}']/\left\langle
f_{1}',\cdots ,f_{m}'\right\rangle,$$

\noindent where $f_{j}'$ is obtained from $f_{j}$ by substituting $g_{i}$ for $y_{i}, 1\leq j\leq m$. \hfill $\square$

\bigskip

The proof of Theorem 4.8 singles out two crucial steps in resolving Weil's problem:

\begin{quote}
i) Find a minimal set of Schubert classes $\left\{ y_{1},\cdots
,y_{n}\right\}$ on $G/P$ that generates the quotient group $
H^{+}(G/P)/D(H^{*}(G/P))$;

ii) Seek a Hilbert basis $\left\{ f_{1},\cdots ,f_{m}\right\}$ of the ideal $\ker \pi$.
\end{quote}

\noindent Since both of the tasks can be implemented by the \textsl{Characteristics} (e.g. \cite[Section 4.4]{DZ7}), we obtain therefore the package \textsl{``Chow-ring''} in MATHEMATICA \cite{DZ4,DZ7} whose function is stated below.

\bigskip

\noindent \textbf{Algorithm III: }\textsl{Chow-ring.}

\begin{quote}
\textbf{Input:} \textsl{The Cartan matrix $C=(a_{{ij}})_{{n\times n}}$ of $G$, and a subset $K\subset \{1,\ldots ,n\}$ to specify a parabolic subgroup $P$}.

\textbf{Output: }\textsl{A presentation (4.4) of the cohomology $H^{\ast
}(G/P)$}.\hfill $\square $
\end{quote}

\bigskip

\noindent \textbf{Example 4.9.} If $G$ is a simple Lie group of rank $n$ and if $K=\{1,\cdots,n\}$, then the parabolic subgroup $P_{K}$ is a maximal torus $T$ on $G$, and the flag manifold $G/T$ is called \textsl{the complete flag manifold} of the group $G$. As applications of the \textsl{Chow-ring} the cohomologies $H^{\ast }(G/T)$ for the exceptional Lie groups have been determined in terms of a minimal system of generators and relations in the Schubert classes (e.g.\cite{DZ6}). We present below the results for $G=F_{4},E_{6},E_{7}$.

\begin{enumerate}
\item[(4.5)] $H^{\ast }(F_{4}/T)=\mathbb{Z}[\omega _{1},\cdots ,\omega
_{4},y_{3},y_{4}]/\left\langle \rho _{2},\rho
_{4},r_{3},r_{6},r_{8},r_{12}\right\rangle $\textsl{, where}

$\rho _{2}=c_{2}-4\omega _{1}^{2}$;

$\rho _{4}=3y_{4}+2\omega _{1}y_{3}-c_{4}$;

$r_{3}=2y_{3}-\omega _{1}^{3}$;

$r_{6}=y_{3}^{2}+2c_{6}-3\omega _{1}^{2}y_{4}$;

$r_{8}=3y_{4}^{2}-\omega _{1}^{2}c_{6}$;

$r_{12}=y_{4}^{3}-c_{6}^{2}$.

\item[(4.6)] $H^{\ast }(E_{6}/T)=\mathbb{Z}[\omega _{1},\ldots ,\omega
_{6},y_{3},y_{4}]/\left\langle \rho _{2},\rho _{3},\rho _{4},\rho
_{5},r_{6},r_{8},r_{9},r_{12}\right\rangle $\textsl{, where}

$\rho _{2}=4\omega _{2}^{2}-c_{2}$;

$\rho _{3}=2y_{3}+2\omega _{2}^{3}-c_{3}$;

$\rho _{4}=3y_{4}+\omega _{2}^{4}-c_{4}$;

$\rho _{5}=2\omega _{2}^{2}y_{3}-\omega _{2}c_{4}+c_{5}$;

$r_{6}=y_{3}^{2}-\omega _{2}c_{5}+2c_{6}$;

$r_{8}=3y_{4}^{2}-2c_{5}y_{3}-\omega _{2}^{2}c_{6}+\omega _{2}^{3}c_{5};$

$r_{9}=2y_{3}c_{6}-\omega _{2}^{3}c_{6}$;

$r_{12}=y_{4}^{3}-c_{6}^{2}$.

\item[(4.7)] $H^{\ast }(E_{7}/T)=\mathbb{Z}[\omega _{1},\ldots ,\omega
_{7},y_{3},y_{4},y_{5},y_{9}]/\left\langle \rho _{i},r_{j}\right\rangle $%
\textsl{, where}

$\rho _{2}=4\omega _{2}^{2}-c_{2}$;

$\rho _{3}=2y_{3}+2\omega _{2}^{3}-c_{3}$;

$\rho _{4}=3y_{4}+\omega _{2}^{4}-c_{4}$;

$\rho _{5}=2y_{5}-2\omega _{2}^{2}y_{3}+\omega _{2}c_{4}-c_{5}$;

$r_{6}=y_{3}^{2}-\omega _{2}c_{5}+2c_{6}$;

$r_{8}=3y_{4}^{2}+2y_{3}y_{5}-2y_{3}c_{5}+2\omega _{2}c_{7}-\omega
_{2}^{2}c_{6}+\omega _{2}^{3}c_{5}$;

$r_{9}=2{y_{9}}+2{y_{4}y_{5}}-2{y_{3}c_{6}}-{\omega _{2}^{2}c_{7}}+{\omega
_{2}^{3}c_{6}}$;

$r_{10}=y_{5}^{2}-2y_{3}c_{7}+\omega _{2}^{3}c_{7}$;

$r_{12}=y_{4}^{3}-4y_{5}c_{7}-c_{6}^{2}-2y_{3}y_{9}-2y_{3}y_{4}y_{5}+2\omega
_{2}y_{5}c_{6}+3\omega _{2}y_{4}c_{7}+c_{5}c_{7}$;

$r_{14}=c_{7}^{2}-2y_{5}y_{9}+2y_{3}y_{4}c_{7}-\omega _{2}^{3}y_{4}c_{7}$;

$%
r_{18}=y_{9}^{2}+2y_{5}c_{6}c_{7}-y_{4}c_{7}^{2}-2y_{4}y_{5}y_{9}+2y_{3}y_{5}^{3}-5\omega _{2}y_{5}^{2}c_{7}
$,
\end{enumerate}

\noindent where the set $\{\omega_{i}, 1\leq i\leq rank(G)\}$ is the Schubert's basis of $H^{2 }(G/T)$, which is also the set of fundamental dominant weights of the relevant group $G$ (e.g.\cite[Lemma 2.4]{D6}); the $c_{r}$'s are certain polynomials in $\omega_{1},\cdots,\omega_{n}$ defined in \cite[(5.17)]{DZ7} whose geometric implication will be made transparent in Example 4.12 below, and where in terms of Bourbaki's numbering of simple roots \cite{Bu}, the $y_{i}$'s are the Schubert classes $s_{I}$ on $G/T$ specified in the table below:

\begin{center}
\begin{tabular}{l|l|l|l|l}
\hline
$y_{i}$ & $y_{3}$ & $y_{4}$ & $y_{5}$ & $y_{9}$ \\ \hline\hline
$F_{4}/T$ & $s_{\{3,2,1\}}$ & $s_{\{4,3,2,1\}}$ &  &  \\ \hline
$E_{6}/T$ & $s_{\{5,4,2\}}$ & $s_{\{6,5,4,2\}}$ &  &  \\ \hline
$E_{7}/T$ & $s_{\{5,4,2\}}$ & $s_{\{6,5,4,2\}}$ & $s_{\{7,6,5,4,2\}}$ & $s_{\{1,5,4,3,7,6,5,4,2\}}$ \\ \hline
\end{tabular}
\end{center}

\bigskip

For more examples of the applications of \textsl{Chow-ring} to computing
with partial flag manifolds $G/P$ we refer to \cite[Theorems 1-7]{DZ4}.\hfill
$\square $
\subsection{Schubert polynomials}
In the classical enumerative geometry, the Chern class $c_{i}$ of the Grassmannian $G_{{n,k}}$ arises firstly as the Poincar\'{e} dual of the variety of the $k$-planes meeting a general $(n-k-i)$-plane on the $n$-space $\mathbb{C}^{n}$, well-known as the $i^{th}$ \textsl{special Schubert class} of $G_{{n,k}}$. The celebrated Giambelli formula \cite[1902]{Gi}, expressing an arbitrary Schubert class $s_{w}$ on $G_{{n,k}}$ as a determinant (i.e. a polynomial) in the special ones, can be praised as the beginning of the idea of Schubert polynomials.

In general, suppose that $G/P$ is a flag manifold for which a solution (4.4) to Weil's problem has been available. Then, every Schubert class $s_{w}$ can be expressed as a polynomial in the generators ${y_{1},\cdots ,y_{n}}$. Inspired by the Giambelli formula, we may call the generators $y_{{1}},\ldots,y_{{n}}$ the \textsl{special Schubert classes} of $G/P$, and ask the question of expressing an arbitrary Schubert class $s_{{w}}$ on $G/P$ as a polynomial $\mathcal{G}_{{w}}$ in the special ones, where $\deg \mathcal{G}_{{w}}=2l(w)$. In particular, starting from Borel's presentation of the cohomology ring of the flag manifold $U(n)/T$ \cite{Bor}

$$H^{*}(U(n)/T)=\mathbb{Z}[x_{1},\cdots,x_{n}]/\langle e_{1},\cdots,e_{n}\rangle,$$

\noindent where $e_{1},\cdots,e_{n}$ are the elementary symmetric polynomials in $x_{1},\cdots,x_{n}$, Lascoux and Sch\"{u}etzenberger \cite{LS} proposed a solution to the problem by the following rules:

\begin{quote}
i) Take the generator $\mathcal{G}_{0}=x_{1}^{n-1}x_{2}^{n-2}\cdots x_{n-1}\in H^{n(n-1)}(U(n)/T)=\mathbb{Z}$ as the top degree Schubert polynomial.

ii) Define $\mathcal{G}_{{w}}:=\mathcal{D}_{ww_{0}}(\mathcal{G}_{0})$, where $w_{0}$ denotes the longest element of the Weyl group $W(U(n))$, and where $\mathcal{D}_{u}$ is the divided difference operator associated to $u\in W(U(n))$ \cite{BGG}.
\end{quote}

\noindent Strictly speaking, Borel's generators $x_{1},\cdots,x_{n-1}$ are the simple roots of the group $SU(n)$, rather than the fundamental dominant weights (i.e. the special Schubert classes of $U(n)/T$). However, this deficiency is not serious, because the linear transformation between these two sets of generators of the ring $H^{*}(U(n)/T)$ is given by the Cartan matrix of $U(n)$ \cite[Lemma 2.3]{D6}. Extending the work of Lascoux-Sch\"{u}etzenberger, the Schubert polynomials of the complete flag manifolds $G/T$ have been defined by Billey-Haiman \cite{BH} for $G=Spin(n),Sp(n)$, and independently by Fomin-Kirillov \cite{FK} for $G=Spin(2n+1).$  In addition, as early as in 1974 Marlin \cite{Ma} has determined the ring $H^{*}(G/T)$ for $G=Spin(n)$ in the context of Schubert calculus.

In the traditional approach to Schubert polynomials, there are two necessary prerequisites (e.g. \cite{BH}):
\begin{quote}
a) Identify the ring $H^{*}(G/P)$ explicitly in terms of a system of special Schubert classes $y_{1},\cdots,y_{n}$ on $G/P$;

b) Specify a polynomial $\mathcal{G}_{0}$ in $y_{1},\cdots,y_{n}$ representing the top degree Schubert class of $G/P$.
\end{quote}
\noindent Nevertheless, granted with \textsl{Characteristics}, we have alternatively a linear algorithm to achieve Schubert polynomials, without resorting to a top degree Schubert polynomial $\mathcal{G}_{0}$, and to the complicated operators $\mathcal{D}_{w}$.

\begin{quote}
\textbf{Algorithm IV:} \textsl{Schubert polynomials}

\textbf{Input:} \textsl{A set} $\{y_{{1}},\ldots ,y_{{n}}\}$ \textsl{of
special Schubert classes on} $G/P$ \textsl{and an integer} $m>0$\textsl{.}

\textbf{Output:} \textsl{a Schubert polynomials} $\mathcal{G}_{{w}}$ \textsl{for} $w\in W^{m}(H;G)$. \hfill $\square $
\end{quote}

We clarify the details of Algorithm IV. Let $\mathbb{Z}[y_{{1}},\ldots ,y_{{n}}]^{(m)}$ be the group of the polynomials of degree $m$ in the special Schubert classes $y_{{1}},\ldots ,y_{{n}}$, and examine the map
\begin{center}
$\pi_{m}:\mathbb{Z}[y_{{1}},\ldots ,y_{{n}}]^{(m)}\rightarrow H^{2m}(G/P)$
\end{center}
\noindent induced by the inclusion $y_{{1}},\ldots ,y_{{n}}\in H^{*}(G/T)$. Let $B(m):=\{y^{\alpha _{{1}}},\ldots ,y^{\alpha _{{b(m)}}}\}$ be the monomial basis of the group $\mathbb{Z}[y_{{1}},\ldots,y_{{n}}]^{(m)}$, and recall by (4.1) that the Schubert basis of the group $H^{2m}(G/P)$ is $\{s_{m,1},\cdots,s_{m,\beta(m)}\}$. Since each $y^{\alpha}\in B(m)$ is a monomial in the special Schubert classes, Algorithm II is applicable to expand it linearly in $\{s_{m,1},\cdots,s_{m,\beta(m)}\}$ to get a $b(m)\times \beta(m)$ matrix $M(\pi_{{m}})$ that satisfies the linear system

\begin{center}
$\left(
\begin{array}{c}
y^{\alpha _{{1}}} \\
\vdots \\
y^{\alpha _{{b(m)}}}%
\end{array}%
\right) =M(\pi _{{m}})\left(
\begin{array}{c}
s_{{m,1}} \\
\vdots \\
s_{{m,\beta (m)}}%
\end{array}%
\right) $.
\end{center}

\noindent Moreover, since the map $\pi_{{m}}$ surjects, the matrix $M(\pi_{{m}})$ has a
$\beta (m)\times \beta (m)$ minor equal to $\pm 1$. Thus, the standard
integral row and column operation diagonalizing $M(\pi _{{m}})$
(\cite[p.162-164]{Sch}) provides us with two invertible matrices
$P=P_{{b(m)\times b(m)}}$ and $Q=Q_{{\beta (m)\times \beta (m)}}$ satisfying the relation

\begin{enumerate}
\item[(4.8)] $PM(\pi _{{m}})Q=\left(
\begin{array}{c}
I_{{\beta (m)}} \\
C%
\end{array}%
\right) _{{b(m)\times \beta (m)}},$
\end{enumerate}

\noindent where $I_{{\beta (m)}}$ denotes the identity matrix of rank $\beta (m)$. Summarizing, Algorithm IV can be realized by following procedure.

\begin{quote}
\textbf{Step 1.} Compute $M(\pi _{{m}})$ using the \textsl{Characteristics};

\textbf{Step 2.} Diagonalize $M(\pi_{{m}})$ to get the matrices $P$ and $Q$ in (4.8);

\textbf{Step 3.} Set $\left(
\begin{array}{c}
\mathcal{G}_{{m,1}} \\
\vdots \\
\mathcal{G}_{{m,\beta(m)}}%
\end{array}%
\right): =Q\cdot \lbrack P]\left(
\begin{array}{c}
y^{\alpha _{{1}}} \\
\vdots \\
y^{\alpha _{{b(m)}}}%
\end{array}%
\right) $,
\end{quote}

\noindent where $[P]$ is the $\beta(m)\times \beta(m)$ matrix formed by the first $\beta(m)$ rows of $P$. It is clear that

\bigskip

\noindent \textbf{Theorem 4.11.} \textsl{The polynomial} $\mathcal{G}_{m,k}(y_{1},\ldots,y_{n})$ \textsl{is a Schubert polynomial of the Schubert class} $s_{m,k}, 1\leq k\leq \beta(m)$. \hfill $\square$

\bigskip

\noindent \textbf{Example 4.12.} Algorithm IV is suitable for those flag manifolds $G/P$ whose top degree Schubert polynomials are absent (or in question).

For the exceptional Lie group $G=E_{n}$ with $n=6,7,8$ the parabolic subgroup $P_{\{2\}}$ has a canonical $n$-dimensional complex representation, that gives rise to the canonical complex $n$-bundle $\xi_{n}$ on the flag manifold $E_{n}/P_{\{2\}}$ \cite{AH}. According to Borel-Hirzebruch \cite[Section 10]{BoH} the Chern classes $c_{i}(\xi_{n})$ can be expressed as a (rather lengthy) polynomial in positive roots of the group $E_{n}$. However, with respect to the special Schubert classes of $E_{n}/P_{\{2\}}$ specified in the table

\begin{center}
\begin{tabular}{l|l}
\hline
$y_{i}$ & $E_{n}/P_{\{2\}},\text{ }n=6,7,8$ \\ \hline\hline
$y_{1}$ & $s_{\{2\}}\text{, }n=6,7,8$ \\ \hline
$y_{3}$ & $s_{\{5,4,2\}}\text{, }n=6,7,8$ \\ \hline
$y_{4}$ & $s_{\{6,5,4,2\}}\text{, }n=6,7,8$ \\ \hline
$y_{5}$ & $s_{\{{7,6,5,4,2}\}}\text{, }n=7,8$ \\ \hline
$y_{6}$ & $s_{\{1,3,6,5,4,2\}}\text{, }n=6,7,8$ \\ \hline
$y_{7}$ & $s_{\{1,3,7,6,5,4,2\}}$, $n=7,8$ \\ \hline
$y_{8}$ & $s_{\{1,3,8,7,6,5,4,2\}}$, $n=8$ \\ \hline
\end{tabular}
\end{center}

\noindent we get from Algorithm IV the following concise expressions of the Chern classes $c_{i}(\xi_{n})$ as polynomials in the special Schubert classes:

\begin{center}
{\footnotesize  \renewcommand{\arraystretch}{0.8}
\begin{tabular}{l|l|l|l}
\hline
& $E_{6}/P_{\{2\}}$ & $E_{7}/P_{\{2\}}$ & $E_{8}/P_{\{2\}}$ \\ \hline\hline
$c_{1}$ & $3y_{1}$ & $3y_{1}$ & $3y_{1}$ \\ \hline
$c_{2}$ & $4y_{1}^{2}$ & $4y_{1}^{2}$ & $4y_{1}^{2}$ \\ \hline
$c_{3}$ & $2y_{3}+2y_{1}^{3}$ & $2y_{3}+2y_{1}^{3}$ & $2y_{3}+2y_{1}^{3}$ \\
\hline
$c_{4}$ & $3y_{4}+y_{1}^{4}$ & $3y_{4}+y_{1}^{4}$ & $3y_{4}+y_{1}^{4}$ \\
\hline
$c_{5}$ & $3y_{1}y_{4}-2y_{1}^{2}y_{3}+y_{1}^{5}$ & $
2y_{5}+3y_{1}y_{4}-2y_{1}^{2}y_{3}+y_{1}^{5}$ & $2y_{5}+3y_{1}y_{4}
-2y_{1}^{2}y_{3}+y_{1}^{5}$ \\ \hline
$c_{6}$ & $y_{6}$ & $y_{6}+2y_{1}y_{5}$ & $5y_{6}+2y_{3}^{2}+6y_{1}
y_{5}-6y_{1}^{2}y_{4}$\\
& & & \ $+4y_{1}^{3}y_{3}-2y_{1}^{6}$ \\ \hline
$c_{7}$ & $0$ & $y_{7}$ & $y_{7}+4y_{1}y_{6}+2y_{1}y_{3}^{2}
+4y_{1}^{2}y_{5}$\\
& & & $\ -6y_{1}^{3}y_{4}+4y_{1}^{4}y_{3}-2y_{1}^{7}$ \\ \hline
$c_{8}$ & $0$  & $0$  & $y_{8}$ \\ \hline
\end{tabular}
}
\end{center}

\noindent By the way, the Schubert polynomials of the flag manifold $E_{6}/P_{\{2\}}$ in degrees $m=8,9$ are

\begin{center}
{\footnotesize  \renewcommand{\arraystretch}{0.8}
\begin{tabular}{l|l}
\hline
$s_{8,1}=y_{4}^{2}-2y_{4}y_{3}y_{1}+y_{4}y_{1}^{4}$ & $%
s_{9,1}=-y_{6}y_{3}+2y_{4}^{2}y_{1}-2y_{4}y_{3}y_{1}^{2}+y_{4}y_{1}^{5}$ \\
\hline
$s_{8,2}=$\ $y_{4}^{2}$ & $s_{9,2}=y_{6}y_{3}-y_{4}^{2}y_{1}$ \\ \hline
$%
s_{8,3}=2y_{4}^{2}-3y_{4}y_{3}y_{1}-y_{4}y_{1}^{4}+3y_{3}^{2}y_{1}^{2}-y_{3}y_{1}^{5}
$ & $s_{9,3}=y_{6}y_{3}-y_{4}^{2}y_{1}+2y_{4}y_{3}y_{1}^{2}-y_{4}y_{1}^{5}$
\\ \hline
$s_{8,4}=2y_{4}^{2}-5y_{4}y_{3}y_{1}+5y_{3}^{2}y_{1}^{2}-2y_{3}y_{1}^{5}$ & $%
s_{9,4}=-y_{6}y_{3}-y_{4}y_{1}^{5}+y_{3}^{3}$ \\ \hline
$%
s_{8,5}=-5y_{4}^{2}+8y_{4}y_{3}y_{1}-y_{4}y_{1}^{4}-5y_{3}^{2}y_{1}^{2}+2y_{3}y_{1}^{5}
$ & $s_{9,5}=-y_{6}y_{3}-3y_{4}^{2}y_{1}+3y_{4}y_{3}y_{1}^{2}-y_{4}y_{1}^{5}$
\\ \hline
\end{tabular}
}
\end{center} \hfill $\square$

\bigskip

\noindent \textbf{Remark 4.13.} Presumably, Schubert polynomials may be useful to compute the characteristics numbers. However, in our approach the Schubert's characteristics is a preliminary step toward Schubert polynomials (e.g. Algorithm IV).

Currently, the theory of Schubert polynomials is a powerful tool for discovering combinatorial structures of the Littlewood-Richardson coefficients \cite{Cos, Buc1,Buc2}, and is essential to the geometric topic of degeneracy loci of maps between vector bundles \cite{FP,Ta}. There have also been extensive studies of Schubert polynomials in quantum cohomology \cite{FGP}, and in the $K$-theory of flag manifolds. For recent progresses of this branch of contemporary Schubert calculus, we refer to the articles Kirillov-Narus \cite{KN}, and Smirnov-Tutubalina \cite{ST}. \hfill $\square$
\subsection{Applications to the topology of homogeneous spaces}
For a compact Lie group $G$ with a closed subgroup $H$ the quotient space $G/H$ is called a \textsl{homogeneous space} of $G$. In contrast to the flag manifolds the cohomology of a homogeneous space may be nontrivial in odd degrees, and may contain torsion elements.

A classical problem of topology is to express the cohomology of a Lie group $G$, or a homogeneous
space $G/H$, by a minimal system of explicit generators and relations. The
traditional approaches due to H.Cartan, A.Borel, P.Baum and H.Toda utilize
various spectral sequence techniques \cite{Bau, Bor, HMS, T, Wo}, and the
calculation encounters the same difficulties when applied to a Lie group $G$
whose integral cohomology has torsion elements, in particular when $G$ is one of the
exceptional Lie groups.

Schubert calculus makes the cohomology theory of homogeneous spaces appearing in a new light. For examples, inputting the formulae $(4.5)-(4.7)$ into the second page of the Serre spectral sequence of the fibration $G\rightarrow G/T$

\begin{equation*}
E_{2}^{**}(G)=H^{*}(G/T)\otimes H^{*}(T),
\end{equation*}

\noindent the integral cohomology $H^{*}(G)$, as well as the Hopf algebra
structure on the $\mod p$ cohomology $H^{*}(G;\mathbb{Z}_{p})$, has been determined
by computing with the Schubert classes on $G/T$ \cite{DZ3,DZ8}. For more examples of the extension of Schubert calculus to computing with the homogeneous spaces, see in \cite[Sect.5]{DZ4}.
\section{Concluding remarks}

Throughout the ages a common hope of geometers is to find calculable
mechanisms among the geometric entities they are caring of (e.g. algebraic
varieties, cellular complexes, vector bundles, or the cobordism classes of
smooth manifolds). The emergence of Schubert's calculus, or the birth of the
intersection theory, had catered to this demand. Today, Schubert's calculus
has widely integrated into many branches of mathematics, and has profoundly affected the trajectories of the development of those fields, such as the theory of characteristic classes \cite{MS}, the string theory \cite{K}, and algebraic combinatorics \cite{Ful0}. All of these vigorously witnessed Hilbert's broad vision and foresight, and at the same time, put forward the essential request to explore effective rules performing the computation.

Subject to the plan this article recalled the earlier studies on Schubert calculus, presented a resolution of the characteristics, and illustrated a passage from the Cartan matrices of Lie groups to the intersection theory of flag manifolds, in which the characteristics play a central role. For the historic significance and rigorous treatments of the enumerative examples of Schubert mentioned in Section 3, we refer to the survey articles S. Kleiman \cite{Kl1,Kl3}, or the relevant sections of the books Fulton \cite{Ful} and Eisenbud-Harris \cite{EH} on intersection theory. For other computer systems that may be used to perform certain computations in the intersection rings of flag manifolds, see for example Nikolenko-Semenov \cite {NS} (the package ChowMaple06),  Grayson et al. \cite {GSS,GS} (the package Schubert2 in Macaulay2), and Decker et al. \cite {DGPS} (the library Schubert in SINGULAR).

\bigskip

\noindent ACKNOWLEDGEMENT. The authors would like to thank their referees for improvements over the earlier version of this paper.

This work is supported by National Science foundation of China: 11771427; 11961131004.

\bigskip
\noindent Haibao Duan, dhb@math.ac.cn

Yau Mathematical Science Center, Tsinghua University, Beijing 100084;

Academy of Mathematics and Systems Sciences, Chinese Academy of Sciences, Beijing 100190.

School of Mathematical Sciences, Dalian University of Technology, Dalian 116024.

\bigskip
\noindent Xuezhi Zhao, zhaoxve@mail.cnu.edu.cn

Department of Mathematics, Capital Normal University, Beijing 100048.

\end{document}